\documentclass{article}
\usepackage{latexsym}
\usepackage{amssymb}
\usepackage{amsmath}
\usepackage{layout}
\usepackage{srcltx}
\newtheorem{teo}{Theorem}[section]
\newtheorem{lema}[teo]{Lemma}

\newtheorem{obs}[teo]{Remark}
\newtheorem{defin}{Definition}
\newenvironment{dem}{\text\bf Proof:}{}
\newcommand{\R}{\mathbb{R}}
\newcommand{\N}{\mathbb{N}}

\newcommand{\deltast}{\Delta_{s,t}}
\newcommand{\xjl}{X^{j,\ell}}
\newcommand{\gamajl}{\gamma^{j,\ell}}
\newcommand{\rojl}{\rho^{j,\ell}}
\newcommand{\xjln}{X^{j,\ell}_n}
\def\C{\mathcal{C}}
\def\bC{\mathbb{C}}

\def\R{\mathbb{R}}

\begin{document}

\begin{center}
\huge
{

Weak convergence to the fractional Brownian sheet from a L\'evy
sheet}

\vspace{0.52cm}

\smallskip

\normalsize {\bf Xavier Bardina\footnote{{$^{,2}$ Partially
supported by Grant PGC2018-097848-B-I00 from MINECO.}} and Carles
Rovira\footnote{ Corresponding author.

{\it E-mail addresses}: Xavier.Bardina@uab.cat;
carles.rovira@ub.edu\hfill}}

\bigskip

{\footnotesize \it $^1$Departament de Matem\`atiques, Universitat
Aut\`onoma de Barcelona, 08193-Bellaterra (Barcelona), Spain.

 $^2$Departament de Matem\`atiques i Inform\`atica, Universitat de Barcelona,
Gran Via 585, 08007-Barcelona, Spain.}

 \end{center}

\begin{abstract}
	In this paper, we show an approximation in law, in the space of the
	continuous functions on $[0,1]^2$, of two-parameter
	Gaussian  processes that can be represented as a
	Wiener type integral by processes constructed from processes that converge to
the Brownian sheet. As an application, we obtain
a sequence of processes constructed from a L\'evy sheet that
converges in law towards the fractional Brownian sheet.

\end{abstract}



\section {Introduction}
%
Let us consider a family of random kernels $\theta_n$ such
that the processes
$$\zeta_n(s,t)=\int_0^t\int_0^s \theta_n(x,y) dxdy,\quad (s,t)\in [0,1]\times [0,1],$$
converge in law in the space of continuous functions $\C([0,1]^2)$
to the Brownian sheet. Our aim is to give sufficient  conditions
on the family $\theta_n$ and on a couple of deterministic kernels $K_1$ and $K_2$ to ensure that the
processes
\begin{equation}\label{prpr}X_{n }(s,t)=\int _{0}^{1}\int
_{0}^{1}K_{1 }(s,u)
K_{2
}(t,v)\theta_n(u,v) dudv,
\end{equation}
converge in law in the space of continuous functions $\C([0,1]^2)$ to the proces  
\begin{equation}
\label{fbs}
W^{K_1,K_2}_{s,t}:= \int _{0}^{1}\int _{0}^{1}
K_{1}(s,u)K_{2}(t,v)dB_{u,v},
\end{equation}
where $B=\left\{B_{s,t},\,(s,t)\in[0,1]^2\right\}$ is a standard
Brownian sheet.

As an example, we obtain the convergence to the fractional Brownian sheet of a family defined using random kernels based on a L\'evy sheet.

 The proces $W^{K_1,K_2}=\left\{W_{s,t}^{K_1,K_2},\, (s,t)\in[0,1]^2
 \right\}$ given by (\ref{fbs})
is characterized by the fact that is centered, Gaussian and its
covariance function factorizes in the following way:
$$E\left[W_{s_1,s_2}^{K_1,K_2}W_{s_1',s_2'}^{K_1,K_2}\right]=\prod_{i=1}^2\left(\int_0^1K_i(s_i,u)K_i(s_i',u)du\right).$$

In the literature there are several paper dealing with the weak convergence to the fractional Brownian motion. In \cite{[DJ]}, \cite{[LD]} and \cite{[SSW]} the authors built up the approximations using Poisson processes while in \cite{[WSSY]} the approximation sequence is based on a  L\'evy process

On the other hand, in \cite{BJ2000}, it is proved  that the family of processes
$$  n\int _{0}^{t}\int _{0}^{s} \sqrt {xy} (-1) ^{N
(\sqrt{n}{x}, \sqrt{n}{y} )}dxdy, \quad n\in\N,
$$ where $\{N(x,y),\,\, (x,y)\in\R_+^2\}$ is a standard Poisson
process in the plane, converges in law in the space of continuous
functions on $[0,1]^2$ to an ordinary Brownian sheet. 
Using this result, in \cite{BJT2003}, the authors show that the sequence  
 $$n\int _{0}^{t}\int _{0}^{s}K_{1 }(s,u)
 K_{2
}(t,v) \sqrt {uv} (-1) ^{N (\sqrt{n}{u}, \sqrt{n}{v})} dudv.
$$
converges in law to the pocess $W^{K_1,K_2 }$ defined in (\ref{fbs}).
Actually, it became a particular case of our Theorem \ref{elteo} since the kernels
$\theta^{0}_n(x,y)=n\sqrt{xy}(-1)^{N(\sqrt{n} x,\sqrt{n}y)}$ will  satisfy our hypothesis (see Section \ref{sec2}).

In \cite{BMQ2020} the result of \cite{BJ2000} is generalized.
The authors consider $\{L(x,y); \, x,y\geq 0\}$ a L\'evy sheet with 
 L\'evy exponent $\Psi(\xi):=a(\xi)+ib(\xi)$, $\xi\in \R$. Given
$\theta\in (0,2\pi)$ such that $a(\theta)a(2\theta)\ne 0$ they define, for any $n\in\N$ and
$(s,t)\in[0,1]^2$,
\begin{equation*}
\bar \zeta_n(s,t):=n K \int_0^{t}\int_0^{s}\sqrt{xy}\, \{\cos(\theta
L(\sqrt{n}x,\sqrt{n}y))+i\sin(\theta L(\sqrt{n}x,\sqrt{n}y))\}dxdy,
\end{equation*}
where the constant $K$ is given by
\begin{equation}\label{eq:50}
K=\frac{1}{\sqrt{2}} \frac{a(\theta)^2+b(\theta)^2}{a(\theta)}.
\end{equation}
Then they prove that, as $n$ tends to
infinity, $\bar \zeta_n$ converges in law, in the space of complex-valued
continuous functions $\C([0,1]^2;\bC)$, to a complex
Brownian sheet. That is, the real part and the imaginary part
converge to two independent Brownian sheets. 

In our paper we show that the random kernels presented in \cite{BMQ2020},
that is,
$$\theta^{1}_n(x,y)=nK\sqrt{xy}\cos(\theta {L(\sqrt{n} x,\sqrt{n}y})$$
and
 $$\theta^{2}_n(x,y)=nK\sqrt{xy}\sin(\theta {L(\sqrt{n} x,\sqrt{n}y})$$ satisfy the set of conditions of Section
\ref{sec2}.
Thus, they can also be used to construct approximations
to the fractional Brownian sheet.

Actually, we will present two sets of conditions (H1) and (H1') on the deterministic kernels $K_1$ and $K_2$. (H1) is satisfied for the kernels that can be used to define the fractional Brownian sheet with parameter bigger that $\frac12$ while the kernels used to define the fractional Brownian sheet with parameter smaller or equal than $\frac12$ satisfy only hypothesis (H1'), that are weaker than (H1). Moreover, deterministic kernels that satisfy only (H1') need random kernels $\theta_n$ that satisfy an extra hypothesis (H4), that depends also on the properties of the deterministic kernels.

We have organized the paper as follows: Section \ref{sec2} is
devoted to present the sets of hypothesis for the deterministic kernels $K_1$ and $K_2$ and for the random kernels $\theta_n$. In Section
\ref{sec3} we prove our main result, that under the hypothesis presented in the previous section we can obtain weak convergence.  In Section \ref{sec4} we prove
that the kernels $\theta^1_n$ and $\theta_n^2$ satisfy the set of
hypothesis, and, finally, in Section \ref{sec5} 
we give some examples for which our result applies, pointing out the case of the fractional Brownian sheet.

Along the paper we will consider two probability space. On one hand
we will consider a probability space $(\Omega,\,\mathcal F,\,P)$,
where we have defined the approximating processes and another
probability space $(\tilde \Omega,\,\tilde {\mathcal F},\,\tilde
P),\,$ where we have defined the limit processes. The mathematical
expectation on these probability spaces will be denoted by $E$ and
$\tilde E$, respectively.

The multiplicative constants that appear along the paper are named
with capital letters and the parameters on which depend are
specified. They may vary from an expression to another one.

\section{Hypothesis}\label{sec2}

Since our aim is to study convergence in law to a Gaussian process
$$W_{s,t}^{K_1,K_2}=\int_0^1\int_0^1K_1(s,u)K_2(t,v)dB_{u,v},$$
where B is a standard Brownian sheet and $K_1$ and $K_2$ are  deterministic kernels we need to fix the conditions that will satisfy $K_1$ and $K_2$ and that will allow us to get the convergence to a fractional Brownian sheet. We consider two sets of hypothesis on $K_1$ and $K_2$, that we recall from \cite{BJ2000}:
\begin{itemize}
\item[ ] {\bf (H1)}\begin{enumerate}
\item[(i)] For $i=1,2$, $K_i$ is measurable and $K_i(0,r)=0$ for all $r\in[0,1]$ almost everywhere.
\item[(ii)] For $i=1,2$, there exists an increasing continuous function
$G_i:[0,1]\longrightarrow\R$ and $\alpha_i>1$ such that for all
$0\leq s<s'\leq1$,
$$\int_0^1\left(K_i(s',r)-K_i(s,r)\right)^2dr\leq\left(G_i(s')-G_i(s)\right)^{\alpha_i}.$$
\end{enumerate}\goodbreak
\item[ ] {\bf (H1')}\begin{enumerate}
\item[(i)] For $i=1,2$, $K_i$ is measurable and $K_i(0,r)=0$ for all $r\in[0,1]$ almost everywhere.
\item[(ii')] For $i=1,2$, there exists an increasing continuous function
$G_i:[0,1]\longrightarrow\R$ and $0<\rho_i\leq1$ such that for all
$0\leq s<s'\leq1$,
$$\int_0^1\left(K_i(s',r)-K_i(s,r)\right)^2dr\leq\left(G_i(s')-G_i(s)\right)^{\rho_i}.$$
\item[(iii)] For $i=1,2$, there exist constants $M_i>0$ and $\beta_i>0$
such that for all $0\leq s<s'\leq1$ and $0\leq s_0<s_0'\leq1$,
$$\int_{s_0}^{s'_0}\left(K_i(s',r)-K_i(s,r)\right)^2dr\leq M_i(s'_0-s_0)^{\beta_i}.$$
\end{enumerate}
\end{itemize}

\begin{obs}
Clearly, condition {\rm (ii)} implies condition {\rm (ii')}. Condition {\rm (iii)}
is added in order to obtain tightness under the weak condition
{\rm (ii')} (see Theorem \ref{elteo}).

For instance, the deterministic kernels associated to fractional Brownian sheet with $H> \frac12$ satisfy {\rm (ii)} while when $H\leq \frac12$ the deterministic kernels satisfy  {\rm (ii')} and {\rm (iii)}.
\end{obs}

On the other hand, to built the approximation sequence $X_n$ given by (\ref{prpr}), we
deal with a set of kernels $\theta_n\in L^{\infty }([0,1]^{2})$.  We consider the following hypothesis on these kernels:
\begin{itemize}
\item[ ] {\bf (H2)}  The processes
$$\zeta_n(s,t)=\int_0^t\int_0^s \theta_n(x,y) dxdy,\quad (s,t)\in [0,1]\times [0,1],$$
converge in law, in the space of continuous functions on $[0,1]^2$,
$\mathcal C([0,1]^2)$, as $n$ tends to infinity, to the Brownian
sheet.
\item[] {\bf (H3)}  For any $f,\,\, g\in L^2([0,1])$
 $$ E\left(\int_{[0,1]^2}f(x)g(y)\theta_n(x,y)dx dy\right)^2\le C \int_{[0,1]^2}f(x)^2g(y)^2 dx dy.$$
\end{itemize} 

Finally, under (H1') we need another hypothesis on the kernels $\theta_n$.
Notice that this hypothesis (H4) also will depend on the properties of $K_1$ and $K_2$. We need  first to introduce some notation.
  Given a real function $X$, defined on $\R_+^2$, and
 $(s,t),(s',t')\in\R_+^2$ such that $s\leq s'$ and $t\leq t'$ we
 denote by $\Delta _{s,t}X(s^{\prime }, t^{\prime })$ the increment
 of $X$ over the rectangle $\left((s,t),(s',t')\right]$, that is $$
 \Delta _{s,t}X(s^{\prime }, t^{\prime })= X(s^{\prime }, t^{\prime
 })-X(s^{\prime }, t)-X(s, t^{\prime })+X(s,t).$$
  When we consider an increment of the processes defined in
 (\ref{prpr}) we have that
 \begin{eqnarray*}
 	\Delta_{st} X_n(s',t')&=& \int_{[0,1]^2}
 	(K_1(s',x)-K_1(s,x))(K_2(t',y)-K_2(t,y))\theta_n(x,y) dx dy
 	\\
 	&=&\Delta_{0,0}Y_n(1,1),
 \end{eqnarray*}
 where the process $Y_n$, that depends on $s,t,s',t',K_1$ and $K_2$,
 is defined by
 \begin{equation}\label{pry}Y_n(s_0,t_0):=\int_{[0,s_0]\times[0,t_0]}
 (K_1(s',x)-K_1(s,x))(K_2(t',y)-K_2(t,y))\theta_n(x,y) dx
 dy.\end{equation}
Now, we can state (H4) when (H1') holds: 
 
\begin{itemize}
\item[] {\bf (H4)} Consider the processes $Y_n$ defined in (\ref{pry}). For any
$0<s_0<s_0'<2s_0$, $0<t_0<t_0'<2t_0$ there exists an even number
$m>\frac{4}{\min\{\rho_1,\rho_2\}}$, such that
\begin{equation}\label{calprovar}
E\left[\Delta_{s_0,t_0}Y_n(s_0',t_0')\right]^m\leq
C_{m,M}\left(G_1(s')-G_1(s)\right)^{\frac{m\rho_1}4}
\left(G_2(t')-G_2(t)\right)^{\frac{m\rho_2}4}\left[(s_0'-s_0)(t_0'-t_0)\right]^{m\gamma},
\end{equation} where $\gamma$
is a parameter belonging to the interval $(0,1)$ (that will only
depend on $\beta_1$ and $\beta_2$) and the constant $C_{m,M}$ only
depends on $m$, $M_1$ and $M_2$.
\end{itemize}

\begin{obs} Notice that under {\rm (H1)} or {\rm (H1')} the processes $X_n$ are continuous.
Indeed, for all $0< s\leq s'<1$, $0<t\leq t'<1$ using condition
(ii') on hypothesis {\rm(H1')} we have that,
$$\left|\Delta_{s,t}X_n(s',t')\right|
\leq\|\theta_n\|_{\infty}\left(G_1(s')-G_1(s)\right)^{\frac{\rho_1}2}\left(G_2(t')-G_2(t)\right)^{\frac{\rho_2}2},$$
where $G_1$ and $G_2$ are continuous functions.
\end{obs}

\section{Convergence in law to two-parameter Gaussian
processes}\label{sec3}

In this section we present our main result. It states as follows:

\begin{teo}\label{elteo} Assume one of the following sets of hypothesis:
	\begin{enumerate}
		\item[{\rm (J1)}]  $K_1$ and $K_2$ satisfy {\rm (H1)} and  the kernels $\theta_n$ satisfy {\rm (H2)} and {\rm (H3)} 
		\item[{\rm (J2)}]  $K_1$ and $K_2$ satisfy {\rm (H1')} and  the kernels $\theta_n$ satisfy {\rm (H2)}, {\rm (H3)}  and {\rm (H4)}.
\end{enumerate}
 Then,  the laws  of the processes
 $\{X_n(s,t),\,\,(s,t)\in [0,1]^2\}$ given by (\ref{prpr}) converge weakly to the law of $\{W^{K_1,K_2}_{s,t}, \,\, (s,t)\in [0,1]^2\}$ in $\mathcal C([0,1]^2)$ when $n$ goes to infinity.
\end{teo}

Before the proof, we need to recall a technical
lemma from \cite{BJT2003} (see Lemma 3.2 therein), that will
be useful for our computations.

\begin{lema}\label{ajus}
	Let $Z=\{Z_{u,v};\,(u,v)\in[0,1]^2\}$ be a continuous process.
	Assume that for a fixed even $m\in\N$ and some
	$\delta_1,{\delta_2}\in(0,1)$ there exists a constant $Q>0$ such
	that
	\begin{equation}\nonumber
	E\left(\Delta_{u,v}Z(u',v')\right)^{m}\leq
	Q\,(u'-u)^{m\delta_1}(v'-v)^{m{\delta_2}}
	\end{equation}
	for any $0<u<u'<2u$, $0<v<v'<2v$. Then, there exists a constant
	$C$ that only depends on $m$, $\delta_1$ and $\delta_2$ such that
	$$
	E\left(\Delta_{u,v}Z(u',v')\right)^{m}\leq C 
	Q\,(u'-u)^{m\delta_1}(v'-v)^{m{\delta_2}}
	$$
	for any $0\leq u<u'\leq1$, $0\leq v<v'\leq1$.
\end{lema}

{\text\bf Proof of Theorem \ref{elteo}:} We will proof the convergence in law checking the tightness of the family of laws of the family $\{X_n\}$ and identifying the limit using the convergence of finite dimensional distributions.

\smallskip

 We first prove the tightness. Using the criterion given by Bickel and
 Wichura in \cite{BW1971} and that our processes are null on the axes it suffices to show
 that for some $m\geq2$ there exist two constants, $C>0$ and
 $\eta>1$, and two increasing continuous functions, $G_1$ and
 $G_2$, such that
 \begin{equation}\label{calprovar2}
\sup_{n}E\left[\Delta_{s,t}X_n(s',t')\right]^m\leq
C\left[\left(G_1(s')-G_1(s)\right)
 \left(G_2(t')-G_2(t)\right)\right]^{\eta},
 \end{equation}
 for any $0\le s \le s' \le 1$, $0\le t \le t'\le 1$.

 Under the set of conditions (J1), using condition (ii) of (H1) and (H3) we have that
 \begin{eqnarray*}
E\left[\Delta_{s,t}X_n(s',t')\right]^2&=&E\left(\int_{[0,1]^2}(K_1(s',x)-K_1(s,x))(K_2(t',y)-K_2(t,y))\theta_n(x,y)dxdy\right)^2\\
&\leq&C\int_{[0,1]^2}(K_1(s',x)-K_1(s,x))^2(K_2(t',y)-K_2(t,y))^2dxdy\\
&\leq&C\left(G_1(s')-G_1(s)\right)^{\alpha_1}\left(G_2(t')-G_2(t)\right)^{\alpha_2},
 \end{eqnarray*}
where $\alpha_1$ and $\alpha_2$ are bigger than 1. So, choosing $\eta=\min\{\alpha_1,\alpha_2\}$,
(\ref{calprovar2}) holds.

Under the set of conditions (J2) we will see
that for $m>\frac{4}{\min\{\rho_1,\rho_2\}}$ given in (H4) and for all $0\leq s<s'\leq1$, $0\leq
t<t'\leq1$,
$$\sup_{n}E\left[\Delta_{s,t}X_n(s',t')\right]^m\leq C\left(G_1(s')-G_1(s)\right)^{\frac{m\rho_1}4}
\left(G_2(t')-G_2(t)\right)^{\frac{m\rho_2}4},$$ where $C$ is a
constant that does not depend on $n$. Indeed,
\begin{eqnarray*}
 &&E\left[ \deltast X_n(s',t')\right]^m\\&=& E\left[\int_{[0,1]^2}
  (K_1(s',x)-K_1(s,x))(K_2(t',y)-K_2(t,y))\theta_n(x,y) dx dy
  \right]^m\\
&=&E\left[\Delta_{0,0}Y_n(1,1)\right]^m,
\end{eqnarray*}
where the process $Y_n$ was defined in (\ref{pry}). By Lemma
\ref{ajus} it suffices to check that  for any
$0<s_0<s_0'<2s_0$, $0<t_0<t_0'<2t_0$ 
$$
E\left[\Delta_{s_0,t_0}Y_n(s_0',t_0')\right]^m\leq
L\left[(s_0'-s_0)(t_0'-t_0)\right]^{m\gamma},
$$ where $L=C\left(G_1(s')-G_1(s)\right)^{\frac{m\rho_1}4}
\left(G_2(t')-G_2(t)\right)^{\frac{m\rho_2}4}$ and it is true by
 hypothesis (H4). So (\ref{calprovar2}) holds easily.

\smallskip

   We proceed now with the identification of the limit law. We will prove the convergence of the finite
   dimensional distributions of the processes $X_n$ to those of $W^{K_1,K_2}$.  Fixed $k\in\N$, consider
   $a_1,\ldots, a_k \in \R$ and $(s_1,t_1),\ldots,(s_k,t_k)\in [0,1]^2$. We must see that the random variables
 \begin{equation}\label{var11}\sum_{j=1}^k a_j X_n(s_j,t_j)\end{equation}
converge in law, as $n$ tends to infinity to
 \begin{equation}\label{var12}\sum_{j=1}^k a_j W^{K_1,K_2}(s_j,t_j).\end{equation}
Actually, we will prove the convergence of the characteristic functions.

 Consider, for any $j$, a sequence $\{ \gamajl\}$ of elementary functions converging in $L^2([0,1])$,
 as $\ell$ tends to infinity, to $K_1(s_j,\cdot)$. In the same way, take a sequence
 $\{\rojl\}$ of elementary functions tending in $L^2([0,1])$ to $K_2(t_j,\cdot)$.
Then, we can introduce the random variables
 $$\xjln=\int_{[0,1]^2}\gamajl(x)\rojl(y)\theta_n(x,y)dx dy,$$
and
 $$\xjl=\int_{[0,1]^2}\gamajl(x)\rojl(y)dB_{x,y}.$$

Then, for any $\lambda\in\R$ we can bound the diference between the characteristic functions of (\ref{var11}) and 
(\ref{var12}) by
\begin{eqnarray}
 && \Big|E\left[e^{i\lambda \sum_j  a_j X_n(s_j,t_j)}   \right]- \tilde E\left[ e^{i\lambda
 \sum_j  a_j
W^{K_1,K_2}(s_j,t_j)}\right]\Big|\nonumber\\ &\le
&\big|E\left[e^{i\lambda \sum_j a_j X_n(s_j,t_j)}  -
 e^{i\lambda\sum_j a_j\xjln
}\right]\Big|\nonumber + \Big| E\left[e^{i\lambda\sum_j
a_j\xjln }\right]-
 \tilde E\left[e^{i\lambda\sum_j a_j\xjl }\right]\Big|\\
&&+ \Big| \tilde E\left[e^{i\lambda\sum_j a_j\xjl } -
e^{i\lambda\sum_i a_j W^{K_1,K_2}(s_j,t_j)}\right]\Big|\nonumber
\\ &:=
&S_{1,n,\ell}+S_{2,n,\ell}+S_{3,\ell}.\label{set}
\end{eqnarray}

We study first $S_{1,n,\ell}$. By the mean value
theorem
$$S_{1,n,\ell} \le C \max_j\{E|X_n(s_j,t_j)-\xjln|\}.$$ Each one of the expectations appearing in the last maximum can be
bounded as follows:
 \begin{eqnarray*}
 &&E|X_n(s_j,t_j)-\xjln|\\
 &=& E\Big|\int_{[0,1]^2}\!
 \!\! K_1(s_j,x)K_2(t_j,y) \theta_n(x,y)
 dx dy-\int_{[0,1]^2}\!\!\! \gamajl(x)\rojl(y)\theta_n(x,y) dx dy\Big|\\
 &\le & E\big|\int_{[0,1]^2}K_1(s_j,x)(\rojl(y)-K_2(t_j,y))\theta_n(x,y)dx dy\Big|\\
 &+ & E\Big|\int_{[0,1]^2}\rojl(y)(K_1(s_j,x)-\gamajl(x))\theta_n(x,y)dx dy \Big|\\
 & \le & C\left(\int_{[0,1]}K_1(s_j,x)^2
dx\right) \left(\int_{[0,1]} (\rojl(y)-K_2(t_j,y))^2 dy \right)\\
 &+& C\left( \int_{[0,1]} (\rojl(y))^2
 dy\right)\left( \int_{[0,1]}(K_1(s_j,x)-\gamajl(x))^2 dx\right),
 \end{eqnarray*}
where in the last inequality we have used hypothesis (H3). Observe that, if $\ell$ is big enough then
 this last expression can be made
arbitrarily small. That is, for any $\varepsilon >0$, there exists $\ell_0$ big enough such that for any $\ell > \ell_0$
\begin{equation}\label{s1}
\sup_{n} | S_{1,n,\ell} | < \varepsilon
\end{equation}

We deal now with $S_{2,n,\ell}$. Since the
functions $\gamajl$ and $\rojl$ are elementary functions, the random
variables $\xjln$ are linear combinations of increments of the
process $\zeta_n(s,t)$  defined in (H2).
The laws of these last processes converge weakly, in $\mathcal
C([0,1]^2)$, to the law of the Brownian sheet, due to (H2). Then,
the linear combinations of the increments of $\zeta_n$ will converge
in law to the same linear combinations of the increments of the
Brownian sheet, that is, to $\xjl$. So, fixed $\ell\in\N$, 
\begin{equation}\label{s2}
\lim_{n \to \infty} S_{2,n,\ell}=0.
\end{equation} 

Finally, we consider $S_{3,\ell}$. Applying the mean value theorem as in the study of $S_{1,n,\ell}$
and using the isometry of the stochastic integral, we can write
\begin{eqnarray*}
S_{3,\ell} & \le & C\max_j\{\tilde E|\xjl-W^{K_1,K_2}(s_j,t_j)|\}\\
 & = & C\max_j\{\tilde E| \int_{[0,1]^2}[\gamajl(x)\rojl(y)-K_1(s_j,x)K_2(t_j,y)]d B_{x,y}|  \} \\
 & \le & C \max_j \left( \int_{[0,1]^2}\big(\gamajl(x)\rojl(y)-K_1(s_j,x)K_2(t_j,y)\big)^2 dx dy \right)^\frac12.
 \end{eqnarray*}
This
last norm in $L^2([0,1]^2)$
 tends to zero, as $\ell$ tends to infinity.
 That is
 \begin{equation}\label{s3}
 \lim_{\ell \to \infty} S_{3,\ell}=0.
 \end{equation}

Putting together (\ref{set}), (\ref{s1}), (\ref{s2}) and (\ref{s3}) we finish the proof.
\hfill $\Box$

\section{Kernels defined from a L\'evy sheet}\label{sec4}

In this section we will prove that the
kernels defined from a L\'evy sheet introduced in \cite{BMQ2020} 
satisfy our hyphotesis.

 \begin{obs}
In \cite{BJ2000} it is proved that if we consider the kernels
 $$\theta^{0}_n(x,y)=n\sqrt{xy}(-1)^{N(\sqrt{n}x,\sqrt{n}y)},$$ where
 $\{N(x,y),\,(x,y)\in\R^2_{+}\}$ is a standard Poisson process then the corresponding
 processes $\{\zeta_n(s,t),\,(s,t)\in[0,1]\times[0,1]\}$ converge in
 law, in the space of continuous funtions on $[0,1]^2$, $\mathcal
 C([0,1]^2)$, to the Brownian sheet. So, hypothesis {\rm (H2)} is verified.
 Moreover, hypothesis {\rm (H3)} corresponds to Lemma 3.1 in
 \cite{BJT2003}, and  {\rm (H4)} for deterministic kernels satisfying {\rm (H1')} is checked in the prove of Lemma
 3.3 in \cite{BJT2003}. Thus, the kernels $\theta_n^0$ satisfy hyphotesis (H2), (H3) and (H4).
 \end{obs}

Let us recall some notation and definitions about L\'evy sheets. If $Q$ is a
rectangle in $\R_+^2$ and $Z$ a random field in $\R_+^2$, we  denote by $\Delta_Q Z$ the increment of $Z$ on $Q$. It is
well-known that, for any negative definite function $\Psi$ in $\R$,
there exists a real-valued random field $L=\{L(s,t); \, s,t\geq 0\}$
such that
\begin{itemize}
\item For any family of disjoint rectangles $Q_1,\dots, Q_n$ in $\R_+^2$, the increments  $\Delta_{Q_1}L,\dots, \Delta_{Q_n}L$ are
independent random variables.
\item For any rectangle $Q$ in $\R_+^2$, the characteristic function of the increment $\Delta_Q L$ is given by
\begin{equation}
\mathbb{E}\left[e^{i \xi \Delta_{Q}L}\right]=
e^{-\lambda(Q)\Psi(\xi)}, \quad \xi\in \R, \label{eq:999}
\end{equation}
where $\lambda$ denotes the Lebesgue measure on $\R_+^2$.
\end{itemize}

\begin{defin}\label{levy}
A random field  $L=\{L(s,t); \, s,t\geq 0\}$  taking values in $\R$
that is continuous in probability and satisfies the above two
conditions is called a L\'evy sheet with exponent $\Psi$.
\end{defin}

By the L\'evy-Khintchine formula, we have
$\Psi(\xi)=a(\xi) + ib(\xi)$, where
$$
a(\xi):=\frac{1}{2}\sigma^2\xi^2 + \int_{\mathbb{R}}[1-\cos(\xi
x)]\eta(dx),
$$
and
$$
b(\xi):=a\xi + \int_{\mathbb{R}}\left[\frac{x \xi}{1+|x|^2} -
\sin(\xi x) \right] \eta(dx).
$$
with $a\in\mathbb{R}$, $\sigma\ge0$ and $\eta$  the corresponding
L\'evy measure, that is a Borel measure on $\R\setminus \{0\}$ that
satisfies
\[
 \int_{\R} \frac{|x|^2}{1+|x|^2} \eta(dx)<\infty.
\]
Notice that $a(\xi)\geq 0$ and, if $\xi\neq 0$, $a(\xi)>0$ whenever
$\sigma>0$ or $\eta$ is nontrivial.

\smallskip
We are able now to recall the kernels defined from a L\'evy sheet introduced in \cite{BMQ2020}. 
Consider $\{L(x,y); \, x,y\geq 0\}$ a L\'evy sheet and
$\Psi(\xi):=a(\xi)+ib(\xi)$, $\xi\in \R$, its L\'evy exponent. Let
$\theta\in (0,2\pi)$  and define
$$\theta^{1}_n(x,y)=nK\sqrt{xy}\cos(\theta {L(\sqrt{n}
x,\sqrt{n}y})$$ and $$\theta^{1}_n(x,y)=nK\sqrt{xy}\sin(\theta
{L(\sqrt{n} x,\sqrt{n}y}),$$ where we assume that
$a(\theta)a(2\theta)\ne 0$ and where the constant $K$ is given by
\begin{equation}
K=\frac{1}{\sqrt{2}} \frac{a(\theta)^2+b(\theta)^2}{a(\theta)}.
\end{equation}

Our aim is to check that these kernels satisfy hyphotesis (H2), (H3) and (H4).
In \cite{BMQ2020} it is proved that the corresponding processes
$\{\zeta_n(s,t),\,(s,t)\in[0,1]\times[0,1]\}$ converge in law, in the
space of continuous funtions on $[0,1]^2$,  to
a Brownian sheet. So, $\theta_n^1$ and $\theta_n^2$ verify
hypothesis (H2).
We will prove that they also satisfy (H3) and (H4) in Lemmas \ref{amel} and \ref{lema2}, respectively. Notice that to check
(H4) we need the additional condition
that $a(\theta)a(2\theta)\cdots a(m\theta)\ne 0$, where  $m$ is the even
appearing in hypothesis (H4). We also present an intermediate technical result in Lemma \ref{lema3}.

\begin{lema}\label{amel}
 For any $f,\,\, g\in L^2([0,1])$ and $k\in{\{1,2\}}$,
 $$ E\left(\int_{[0,1]^2}f(x)g(y)\theta^k_n(x,y)dx dy\right)^2\le \frac{136}{a^2(\theta)} K^2 \int_{[0,1]^2}f(x)^2g(y)^2 dx dy.$$
\end{lema}
\begin{dem}
We will prove the lemma in the case $k=1$ since the case $k=2$ can be done using similar computations. We have that

\begin{eqnarray}
 && E\left(\int_{[0,1]^2}f(x)g(y)\theta^1_n(x,y)dx
dy\right)^2 \nonumber\\ &=& \int_{[0,1]^4}E\left[\prod_{j=1}^2
\left(f(x_j)g(y_j)\theta^1_n(x_j,y_j)\right)\right]dx_1
 \cdots
 dy_2\nonumber\\
&= & 2n^2K^2\int_{[0,1]^4}\prod_{j=1}^2 \left(f(x_j)g(y_j)\sqrt{x_j y_j}\right)\nonumber\\
 &\times & E[\cos(\theta L(\sqrt{n}x_1,\sqrt{n}y_1))\cos(\theta L(\sqrt{n}x_2,\sqrt{n}y_2))]\,
1_{\{x_1\le x_2\}}dx_1\cdots dy_2.\label{bs}
\end{eqnarray}

Notice that, using complex notation, we have that
\begin{eqnarray*}
&&\cos(\theta L(\sqrt{n}x_1,\sqrt{n}y_1))\cos(\theta
L(\sqrt{n}x_2,\sqrt{n}y_2))\\
&=&\left(\frac{e^{i\theta L(\sqrt{n}x_1,\sqrt{n}y_1)}+e^{-i\theta
L(\sqrt{n}x_1,\sqrt{n}y_1)}}{2}\times\frac{e^{i\theta
L(\sqrt{n}x_2,\sqrt{n}y_2)}+e^{-i\theta
L(\sqrt{n}x_2,\sqrt{n}y_2)}}{2}\right)\\
&=&\frac14\left(e^{i\theta L(\sqrt{n}x_1,\sqrt{n}y_1)+i\theta
L(\sqrt{n}x_2,\sqrt{n}y_2)} + e^{i\theta
L(\sqrt{n}x_1,\sqrt{n}y_1)-i\theta
L(\sqrt{n}x_2,\sqrt{n}y_2)}\right.\\&&\left.+e^{-i\theta
L(\sqrt{n}x_1,\sqrt{n}y_1)+i\theta
L(\sqrt{n}x_2,\sqrt{n}y_2)}+e^{-i\theta
L(\sqrt{n}x_1,\sqrt{n}y_1)-i\theta L(\sqrt{n}x_2,\sqrt{n}y_2)}
\right)\\
&:=& A_1+A_2+A_3+A_4.
\end{eqnarray*}
Putting this expression  in (\ref{bs}) we  obtain that the term
(\ref{bs}) is equal to the sum of the corresponding four terms obtained from $A_1, A_2, A_3$ and $A_4$
that we will denote by $I_1$, $I_2$, $I_3$ and $I_4$.

We will deal first with $I_1$ and so, we will begin with the study of $A_1$.
Notice that
\begin{eqnarray*}
&&E\left[e^{i\theta L(\sqrt{n}x_1,\sqrt{n}y_1)+i\theta
L(\sqrt{n}x_2,\sqrt{n}y_2)}I_{\{x_1\leq x_2\}} \right]\\
&=& E\left[e^{i\theta
\Delta_{0,0}L(\sqrt{n}x_1,\sqrt{n}y_1)+i\theta
\Delta_{0,0}L(\sqrt{n}x_2,\sqrt{n}y_2)}I_{\{x_1\leq x_2\}} \right]\\
&=&
e^{-n[(x_2-x_1)y_2+(y_1-y_2)x_1]\Psi(\theta)-n(x_1y_2\Psi(2\theta))}I_{\{x_1\leq
x_2\}}I_{\{y_2\leq y_1\}}\\
&&+e^{-n(x_2y_2-x_1y_1)\Psi(\theta)-n(x_1y_1\Psi(2\theta))}I_{\{x_1\leq
x_2\}}I_{\{y_1\leq y_2\}}\\
&\leq&e^{-n[(x_2-x_1)y_2+(y_1-y_2)x_1]a(\theta)}I_{\{x_1\leq
x_2\}}I_{\{y_2\leq y_1\}}\\
&&+e^{-n[(x_2-x_1)y_1+(y_2-y_1)x_1]a(\theta)}I_{\{x_1\leq
x_2\}}I_{\{y_1\leq y_2\}},
\end{eqnarray*}
where in the last inequality  we have bounded by 1
the modulus of all the terms with the factor $i b(\theta)$ in the
exponential and we have used that $a(\theta) \ge 0$ to bound also by 1 some exponentials with negative real exponent.

Since both summands in the last expression are equal interchanging the roles
of $y_1$ and $y_2$, we obtain that
$$|I_1|\leq n^2K^2\int_{[0,1]^4}\prod_{j=1}^2 \left|f(x_j)g(y_j)\sqrt{x_j y_j}\right|e^{-n[(x_2-x_1)y_1+(y_2-y_1)x_1]a(\theta)}I_{\{x_1\leq
x_2\}}I_{\{y_1\leq y_2\}}dx_1\cdots dy_2.$$

On the other hand, using that $a(-\theta)=a(\theta)$, we can bound
the modulus of the integrals $I_2$, $I_3$ and $I_4$ but the same
bound. Then,
\begin{eqnarray}
 && E\left(\int_{[0,1]^2}f(x)g(y)\theta^1_n(x,y)dx
dy\right)^2\nonumber\\ &\leq& 4n^2K^2\int_{[0,1]^4}\prod_{j=1}^2
\left|f(x_j)g(y_j)\sqrt{x_j
y_j}\right|e^{-n[(x_2-x_1)y_1+(y_2-y_1)x_1]a(\theta)}I_{\{x_1\leq
x_2\}}I_{\{y_1\leq y_2\}}dx_1\cdots
dy_2.\nonumber
\end{eqnarray}

From here we can follow closely the prove of Lemma 3.1 in \cite{BJT2003}. We have to divide
the region of integration in two parts: $A:=\{x_1\leq x_2\leq
2x_1,y_1\leq y_2\leq 2y_1\}$ and $A^C$.
Over the region of integration $A$  we can obtain
the bound
$$2\times\frac4{a^2(\theta)} K^2 \int_{[0,1]^2} \left(f(x)g(y)\right)^2dxdy,$$
while over
 the region $A^C$, 
we get the bound
$$2\times\frac{64}{a^2(\theta)} K^2 \int_{[0,1]^2} \left(f(x)g(y)\right)^2dxdy,$$
and the proof finishes easily.

\hfill$\Box$
\end{dem}

Let us present an intermediate technical result.

\begin{lema}\label{lema3}
	Consider  $L=\{L(s,t); \, s,t\geq 0\}$  
 a L\'evy sheet with exponent 
$\Psi(\xi)=a(\xi) + ib(\xi)$, $0<s_0<s_0'<2s_0$,
$0<t_0<t_0'<2t_0$, $m$ an even number  and $\theta \in (0,2\pi) $ such that $a(\theta)a(2\theta)\cdots a(m\theta)\ne 0$. Then, for any $(x_1,y_1),\ldots, (x_m,y_m)$ such that 
$s_0< x_j <s_0'$ and
$t_0<y_j<t_0'$ for all $j \in \{1,\ldots,m\}$, it holds that
\begin{eqnarray*} & &E\left[ \prod_{j=1}^m\cos(\theta L(\sqrt{n}x_j,\sqrt{n}y_j)) \right] \\
&\leq&\exp\left[-\frac{n}{2}a^{*}(\theta)\big[{x_{(m-1)}}(y_{(m)}-y_{(m-1)})+{x_{(m-3)}}(y_{(m-2)}-y_{(m-3)})+\cdots+
{x_{(1)}}(y_{(2)}-y_{(1)})\big]\right]\\&&\times\exp
\left[-\frac{n}2a^{*}(\theta)\big[{y_{(m-1)}}(x_{(m)}-x_{(m-1)})+{y_{(m-3)}}(x_{(m-2)}-x_{(m-3)})+\cdots+
{y_{(1)}}(x_{(2)}-x_{(1)})\big]\right],
\end{eqnarray*}
where $a^{*}(\theta)=\min\{a(\theta),a(2\theta),\dots,a(m\theta)\}$ and $x_{(1)},\ldots,x_{(m)}$ and $y_{(1)},\ldots,y_{(m)}$ are  $x_1,\ldots,x_m$ and $y_1,\ldots,y_m$ after being ordered.
\end{lema}

\begin{dem}
Notice that,
\begin{eqnarray}
\prod_{j=1}^m\cos(\theta L(\sqrt{n}x_j,\sqrt{n}y_j))
&=&\prod_{j=1}^m\left(\frac{e^{i\theta
L(\sqrt{n}x_j,\sqrt{n}y_j)}+e^{-i\theta
L(\sqrt{n}x_j,\sqrt{n}y_j)}}{2}\right)\nonumber\\
&=&\frac1{2^m}\sum_{(\delta_1,\dots,\delta_m)\in\{1,-1\}^n}\prod_{j=1}^m
e^{i\delta_j\theta L(\sqrt{n}x_j,\sqrt{n}y_j)}\nonumber\\
&=&\frac1{2^m}\sum_{(\delta_1,\dots,\delta_m)\in\{1,-1\}^n}
e^{i\theta \sum_{j=1}^m \delta_j L(\sqrt{n}x_j,\sqrt{n}y_j)}\nonumber\\
&=&\frac1{2^m}\sum_{(\delta_1,\dots,\delta_m)\in\{1,-1\}^n} e^{i\theta
\sum_{j=1}^m \delta_j \Delta_{0,0}L(\sqrt{n}x_j,\sqrt{n}y_j)}.\label{sumdelta}
\end{eqnarray}
Fixed $(\delta_1,\dots,\delta_m)$, we can write
\begin{eqnarray*}
\sum_{j=1}^m \delta_j \Delta_{0,0}L(\sqrt{n}x_j,\sqrt{n}y_j)
&=& \sum_{j=1}^m \delta_j \Delta_{s_0,t_0}L(\sqrt{n}x_j,\sqrt{n}y_j)
+\sum_{j=1}^m \delta_j \Delta_{s_0,0}L(\sqrt{n}x_j,\sqrt{n}t_0)\\&&
+\sum_{j=1}^m \delta_j \Delta_{0,t_0}L(\sqrt{n}s_0,\sqrt{n}y_j)
+L(\sqrt{n}s_0,\sqrt{n}t_0)\sum_{j=1}^m \delta_j.
\end{eqnarray*}
and so
\begin{eqnarray}
e^{\sum_{j=1}^m \delta_j \Delta_{0,0}L(\sqrt{n}x_j,\sqrt{n}y_j)}
&=& e^{\sum_{j=1}^m \delta_j
\Delta_{s_0,t_0}L(\sqrt{n}x_j,\sqrt{n}y_j)}e^{\sum_{j=1}^m \delta_j
\Delta_{s_0,0}L(\sqrt{n}x_j,\sqrt{n}t_0)} \nonumber\\&&\times e^{\sum_{j=1}^m
\delta_j
\Delta_{0,t_0}L(\sqrt{n}s_0,\sqrt{n}y_j)}e^{L(\sqrt{n}s_0,\sqrt{n}t_0)\sum_{j=1}^m
\delta_j}.\label{indepe}
\end{eqnarray}

Since the  interval $((0,0),(s_0,t_0)]$  and the families of intervals
$\{((s_0,t_0),(x_j,y_j)], 1 \le j \le m  \}$, 
$\{((s_0,0),(x_j,t_0)], 1 \le j \le m  \}$,
$\{((0,t_0),(s_0,y_l)], 1 \le j \le m  \}$
 have support on disjoint sets,
the four factors of (\ref{indepe}) are independent random
variables. Moreover,
\begin{eqnarray}
&&|E\left(e^{i\theta\sum_{j=1}^m \delta_j \Delta_{0,0}L(\sqrt{n}x_j,\sqrt{n}y_j)}\right)| \nonumber\\
&=& |E\left(e^{i\theta\sum_{j=1}^m \delta_j
\Delta_{s_0,t_0}L(\sqrt{n}x_j,\sqrt{n}y_j)}\right)||E\left(e^{i\theta\sum_{j=1}^m
\delta_j
\Delta_{s_0,0}L(\sqrt{n}x_j,\sqrt{n}t_0)}\right)|  \nonumber\\&&\times|E\left(e^{i\theta\sum_{j=1}^m
\delta_j
\Delta_{0,t_0}L(\sqrt{n}s_0,\sqrt{n}y_j)}\right)||E\left(e^{i\theta
L(\sqrt{n}s_0,\sqrt{n}t_0)\sum_{j=1}^m
\delta_j}\right)| \nonumber\\
&\leq&|E\left(e^{i\theta\sum_{j=1}^m \delta_j
\Delta_{s_0,0}L(\sqrt{n}x_j,\sqrt{n}t_0)}\right)||E\left(e^{i\theta\sum_{j=1}^m
\delta_j \Delta_{0,t_0}L(\sqrt{n}s_0,\sqrt{n}y_j)}\right)|, \label{twoterms}
\end{eqnarray}
where in the last step we have bounded two factors by 1.

Let us study first the second term in (\ref{twoterms}).
We need to introduce the notation $(y_{(1)},\delta_{(1)}),\dots,(y_{(n)},\delta_{(n)})$ for the
variables $(y_{1},\delta_{1}),\dots,(y_{n},\delta_{n})$ sorted in increasing order by the variables $y_i$. 
Using this notation, we can write
\begin{eqnarray*}
&&\sum_{j=1}^m \delta_j
\Delta_{0,t_0}L(\sqrt{n}s_0,\sqrt{n}y_j) = \sum_{j=1}^m \delta_{(j)}
\Delta_{0,t_0}L(\sqrt{n}s_0,\sqrt{n}y_{(j)})  \\
&=&\delta_{(m)} \Delta_{0,y_{(m-1)}}L(\sqrt{n}s_0,\sqrt{n}y_{(m)})+
(\delta_{(m)}+\delta_{(m-1)})
\Delta_{0,y_{(m-2)}}L(\sqrt{n}s_0,\sqrt{n}y_{(m-1)})\\&&+\dots
+\left(\sum_{j=1}^m \delta_{(j)}\right)
\Delta_{0,t_0}L(\sqrt{n}s_0,\sqrt{n}y_{(1)}).
\end{eqnarray*}
Since in the last expression all the rectangles where we consider the increments of the L\'evy sheet are disjoint,  all the terms in the last expression are independent random variables. Then, if we
bound by 1 all the factors with an even number of summands in
$\delta_{(m)}+\delta_{(m-1)}+\dots+\delta_{(j)}$  we obtain that
\begin{eqnarray*}
&&|E\left(e^{i\theta\sum_{j=1}^m \delta_{j}
\Delta_{0,t_0}L(\sqrt{n}s_0,\sqrt{n}y_{j})}\right)| \\
&\leq&
\exp\big[-ns_0(y_{(m)}-y_{(m-1)})\Psi(\delta_{(m)}\theta)-ns_0(y_{(m-2)}-y_{(m-3)})\Psi((\delta_{(m)}+\delta_{(m-1)}+\delta_{(m-2)})\theta) \\&&-\cdots-
ns_0(y_{(2)}-y_{(1)})\Psi((\sum_{j=2}^m\delta_{(j)})\theta)\big]. 
\end{eqnarray*}
Let us recall that  $\Psi(h\theta)=a(h\theta)+ib(h\theta)$ for all
$h\in\R$. Bounding again by 1 the modulus of all the terms with the factor $i
b(h\theta)$ in the exponential we obtain that 
\begin{eqnarray}&&|E\left(e^{i\theta\sum_{j=1}^m \delta_{j}
		\Delta_{0,t_0}L(\sqrt{n}s_0,\sqrt{n}y_{j})}\right)|\nonumber\\
&\leq&
\exp\big[-ns_0(y_{(m)}-y_{(m-1)})a(\delta_{(m)}\theta)-ns_0(y_{(m-2)}-y_{(m-3)})a((\delta_{(m)}+\delta_{(m-1)}+\delta_{(m-2)})\theta)\nonumber\\&&-\cdots-
ns_0(y_{(2)}-y_{(1)})a((\sum_{j=2}^m\delta_{(j)})\theta)\big]\nonumber\\
&\leq&
\exp\left[-ns_0(y_{(m)}-y_{(m-1)})a^{*}(\theta)-ns_0(y_{(m-2)}-y_{(m-3)})a^{*}(\theta)-\cdots-
ns_0(y_{(2)}-y_{(1)})a^{*}(\theta)\right],\nonumber\\&& \label{termi2}
\end{eqnarray}
where $a^{*}(\theta)=\min\{a(\theta),a(2\theta),\dots,a(m\theta)\}$.
\smallskip

Using the same arguments we can bound also the first term in (\ref{twoterms}) and we get that
\begin{eqnarray}
&&|E\left(e^{i\theta\sum_{j=1}^m \delta_j
\Delta_{s_0,0}L(\sqrt{n}x_j,\sqrt{n}t_0)}\right)|\nonumber\\
&\leq&
\exp\left[-nt_0(x_{(m)}-x_{(m-1)})a^{*}(\theta)-nt_0(x_{(m-2)}-x_{(m-3)})a^{*}(\theta)-\cdots-
nt_0(x_{(2)}-x_{(1)})a^{*}(\theta)\right].\nonumber \\&& \label{termi1}
\end{eqnarray}

Then, putting together (\ref{twoterms}), (\ref{termi1}) and (\ref{termi2})
\begin{eqnarray*}
&&|E\left(e^{i\theta\sum_{j=1}^m \delta_j \Delta_{0,0}L(\sqrt{n}x_j,\sqrt{n}y_j)}\right)|\\
&\leq&\exp\left[-na^{*}(\theta)s_0\big[(y_{(m)}-y_{(m-1)})+(y_{(m-2)}-y_{(m-3)})+\cdots+
(y_{(2)}-y_{(1)})\big]\right]\\&&\times
\exp\left[-na^{*}(\theta)t_0\big[(x_{(m)}-x_{(m-1)})+(x_{(m-2)}-x_{(m-3)})+\cdots+
(x_{(2)}-x_{(1)})\big]\right].
\end{eqnarray*}

Finnally, using that $2t_0>t'_0$ and $2s_0>s'_0$, the last expression can be
bounded by
\begin{eqnarray}
&&\exp\left[-na^{*}(\theta)\frac{s'_0}{2}\big[(y_{(m)}-y_{(m-1)})+(y_{(m-2)}-y_{(m-3)})+\cdots+
(y_{(2)}-y_{(1)})\big]\right]\nonumber\\&&\times
\exp\left[-na^{*}(\theta)\frac{t'_0}{2}\big[(x_{(m)}-x_{(m-1)})+(x_{(m-2)}-x_{(m-3)})+\cdots+
(x_{(2)}-x_{(1)})\big]\right]\nonumber\\
&\leq&\exp\left[-\frac{n}{2}a^{*}(\theta)\big[{x_{(m-1)}}(y_{(m)}-y_{(m-1)})+{x_{(m-3)}}(y_{(m-2)}-y_{(m-3)})+\cdots+
{x_{(1)}}(y_{(2)}-y_{(1)})\big]\right]\nonumber\\&&\times\exp
\left[-\frac{n}2a^{*}(\theta)\big[{y_{(m-1)}}(x_{(m)}-x_{(m-1)})+{y_{(m-3)}}(x_{(m-2)}-x_{(m-3)})+\cdots+
{y_{(1)}}(x_{(2)}-x_{(1)})\big]\right].\nonumber\\&& \label{fini}
\end{eqnarray}
Since this last bound does not depend on $(\delta_1,\ldots,\delta_m)$, putting together the bound (\ref{fini}) with (\ref{sumdelta})
we finish the proof easily.
\hfill$\Box$
\end{dem}

\smallskip

We are able now to proof (H4) under (H1').

\begin{lema}\label{lema2} Let us consider
	the processes (\ref{pry}) defined using the kernels $\theta_n^1$ or $\theta_n^2$, with  $\theta \in (0,2\pi) $ such that $a(\theta)a(2\theta)\cdots a(m\theta)\ne 0$. Assume {\rm (H1')}.  Then, for any $0<s_0<s_0'<2s_0$,
	$0<t_0<t_0'<2t_0$ there exists an even number
	$m>\frac{4}{\min\{\rho_1,\rho_2\}}$ such that
	\begin{equation}\label{calprovar3}
	E\left[\Delta_{s_0,t_0}Y_n(s_0',t_0')\right]^m\leq
	C_{m,M}\left(G_1(s')-G_1(s)\right)^{\frac{m\rho_1}4}
	\left(G_2(t')-G_2(t)\right)^{\frac{m\rho_2}4}\left[(s_0'-s_0)(t_0'-t_0)\right]^{m\gamma},
	\end{equation}  where $\gamma$
	is a parameter belonging to the interval $(0,1)$ (that will only
	depend on $\beta_1$ and $\beta_2$) and the constant $C_{m,M}$ only
	depends on $m$, $M_1$ and $M_2$.
\end{lema}

\begin{dem}
	As in lemma \ref{amel}, we will prove the result only for
	the kernels $\theta_n^1$ since the case $\theta_n^2$ is very similar. From ({\ref{pry}}), the definition of $Y_n$, we can write that
	\begin{eqnarray*}
		&&E\left[\Delta_{s_0,t_0}Y_n(s_0',t_0')\right]^m\\
		&= & n^m K^m
		E\big[\int_{[0,1]^{2m}}\prod_{j=1}^m \Bigg( I_{[s_0',s_0]}(x_j)I_{[t_0',t_0]}(y_j)
		(K_1(s',x_j)-K_1(s,x_j))(K_2(t',y_j)-K_2(t,y_j))\\&&\times\sqrt{x_j
			y_j}\cos(\theta L(\sqrt{n}x_j,\sqrt{n}y_j))  \Bigg)
		dx_1\cdots dy_m   \big].
	\end{eqnarray*}
Using lemma \ref{lema3}, we obatin easily that
\begin{eqnarray*}
	&&E\left[\Delta_{s_0,t_0}Y_n(s_0',t_0')\right]^m\\
	&\leq&(m!)^2n^mK^m
	\int_{[0,1]^{2m}}\prod_{j=1}^m  \Bigg( 
	I_{[s_0',s_0]}(x_j)I_{[t_0',t_0]}(y_j)\\&&\times
	(K_1(s',x_j)-K_1(s,x_j))(K_2(t',y_j)-K_2(t,y_j))\sqrt{x_j y_j}  \Bigg) 
	\\&&\times
	\exp\left[-\frac{n}2a^{*}(\theta)\left[x_{m-1}(y_{m}-y_{m-1})+\cdots+
	x_{1}(y_{2}-y_{1})\right]\right] \\&&\times\exp\left[-\frac{n}2
	a^{*}(\theta)y_{m-1}(x_{m}-x_{m-1})+\cdots+
	y_{1}(x_{2}-x_{1})\right]I_{\{x_1\leq\cdots\leq
		x_m\}}\\&&\times I_{\{y_1\leq\cdots\leq y_m\}}dx_1\cdots dy_m\\
&\leq&C_m \Bigg( n^2K^2
\int_{[0,1]^{4}}\prod_{j=1}^2  \Bigg(  I_{[s_0',s_0]}(x_j)I_{[t_0',t_0]}(y_j)\\&&\times
(K_1(s',x_j)-K_1(s,x_j))(K_2(t',y_j)-K_2(t,y_j))\sqrt{x_j y_j}  \Bigg)
\\&&\times \exp\left[-\frac{n}2a^{*}(\theta)\left[
x_{1}(y_{2}-y_{1})+y_{1}(x_{2}-x_{1})\right]\right]I_{\{x_1\leq
x_2\}}I_{\{y_1\leq y_2\}}dx_1\cdots dy_2 \Bigg)^{\frac{m}2},
\end{eqnarray*}
where $C_m$ is a constant depending only on $m$.

By the computations of Lemma \ref{amel} this last expression is
bounded by

 \begin{eqnarray*}
 & &C_m \frac{K^2}{(a^*(\theta))^2}\Big( \int_{[0,1]^2}I_{[s_0',s_0]}(x)I_{[t_0',t_0]}(y) (K_1(s',x)-K_1(s,x))^2
 (K_2(t',y)-K_2(t,y))^2
   dx dy\Big)^{\frac{m}2}\\
   &\leq&C_m\frac{K^2}{(a^*(\theta))^2}\Big( \int_{[0,1]^2}I_{[s_0',s_0]}(x)I_{[t_0',t_0]}(y) (K_1(s',x)-K_1(s,x))^2
   (K_2(t',y)-K_2(t,y))^2
   dx dy\Big)^{\frac{m}4}\\
   &&\times\Big( \int_{[0,1]^2}(K_1(s',x)-K_1(s,x))^2  (K_2(t',y)-K_2(t,y))^2
   dx dy\Big)^{\frac{m}4}.
\end{eqnarray*}

Using conditions (iii) and (ii') on the kernels $K_1$ and $K_2$ the
last expression is bounded by
\begin{eqnarray*}
&&C_m M_1^{\frac m4} M_2^{\frac m4} \frac{K^2}{(a^*(\theta))^2} (s_0'-s_0)^{\frac{\beta_1
m}4}(t_0'-t_0)^{\frac{\beta_2 m}4}
\left(G_1(s')-G_1(s)\right)^{\frac{m\rho_1}4}
   \left(G_2(t')-G_2(t)\right)^{\frac{m\rho_2}4}\\
   &\leq& C_{m,M} \frac{K^2}{(a^*(\theta))^2} \left[(s_0'-s_0)(t_0'-t_0)\right]^{\gamma m}
   \left(G_1(s')-G_1(s)\right)^{\frac{m\rho_1}4}
   \left(G_2(t')-G_2(t)\right)^{\frac{m\rho_2}4},
   \end{eqnarray*}
where $\gamma=\frac14\inf\{\beta_1,\beta_2\}$. So, we have proved
inequality (\ref{calprovar3}) and the proof is now complete.

\hfill $\Box$
\end{dem}

\section{The fractional Brownian sheet and other examples}\label{sec5}

The main example that satisfies our hypothesis and with which we can
apply Theorem \ref{elteo} is the fractional Brownian sheet. More
precisely, we use the anisotropic fractional Wiener random field,
introduced by \cite{K1996} and \cite{ALP2002}. This is a centered
Gaussian process, defined on some probability space
 $(\tilde \Omega, \tilde{\mathcal F}
, \tilde P )$, denoted by $W^{\alpha,
\beta}=\{W^{\alpha,\beta}_{s,t},\,\, (s,t)\in \R_+^2\}$, with
covariance function
 given by
\begin{equation}
  \label{cov}
  \tilde E \left(  W^{\alpha, \beta }_{s,t}  W^{\alpha, \beta }_{s^{\prime },t^{\prime }}\right)
  =\frac{1}{2}\left( {s'}^{2\alpha} +s^{2\alpha} -|s'-s|^{2\alpha}\right)\frac{1}{2}\left(
  {t'}^{2\beta} +t^{2\beta} -
  {|t'-t|}^{2\beta}\right),
\end{equation}
  where $\alpha $ and $\beta$ are two parameters belonging to the interval $(0,1)$.
By definition, this process is null almost surely on the axes, and
it is proved in \cite{K1996} and \cite{ALP2002} that it possesses a
continuous version. Observe that if $\alpha=\beta=\frac12$, then we
obtain the standard Brownian sheet.

Recall that a fractional Brownian motion of Hurst parameter
$\alpha\in (0,1)$, $\,\,W^{\alpha}=\{W^{\alpha}_t, \,\,
t\in\R_+\}$
 is
 a centered Gaussian process with
covariance function
 $$
 R(t,s)=E\left( W^{\alpha}_{t}W^{\alpha}_{s}\right) =\frac{1}{2}\left( t^{2\alpha} +s^{2\alpha}
 -|t-s|^{2\alpha}\right).
 $$
 This process admits an integral representation of the form (see for instance \cite{AMN2001})
 \begin{equation}\label{eq1}
W_t^{\alpha}=\int_0^t K_{\alpha}(t,s)dB_s,
\end{equation}
where $B$ is a standard Brownian motion and the kernel
$K_{\alpha}$ is given by
\begin{equation}\label{eq2}
K_{\alpha}(t,s)=\left[d_{\alpha} (t-s)^{\alpha-\frac12}+
d_{\alpha}(\frac12-\alpha)\int_s^t
(u-s)^{\alpha-\frac32}\left(1-\left(\frac{s}{u}\right)^{\frac12-\alpha}\right)
du\right]I_{(0,t)}(s),
\end{equation} with $d_{\alpha}$ the following normalizing constant
 $$ d_{\alpha}=\left(\frac{2\alpha\Gamma(\frac32-\alpha)}{\Gamma(\alpha+\frac12)
 \Gamma(2-2\alpha)}\right)^{\frac12}.$$

%

  Taking into account the expression (\ref{eq1}) for the fractional Brownian motion, one can consider
 the following representation in law of the fractional Brownian
 sheet:
\begin{eqnarray*}
   \int _{0}^{t}\int _{0}^{s}
   K_{\alpha }(s,u)K_{\beta }(t,v)dB_{u,v}.
   \end{eqnarray*}
To see the equality in law of this processes and the fractional Brownian sheet it suffices
   to realize that this is a centered Gaussian process with the same covariance function as $W^{\alpha,\beta}$.

If we consider the kernels $K_{\alpha}$ (or $K_{\beta}$) we have
that
$$
\int_0^1\left(K_{\alpha}(s',r)-K_{\alpha}(s,r)\right)^2dr=E\left(W_{s'}^{\alpha}-
W_s^{\alpha}\right)^2 =(s'-s)^{2\alpha}.
$$

And then the set of conditions (H1) is satisfied if $\alpha>\frac12$
and $\beta>\frac12$. In \cite{BJT2003} it is proved that if $\alpha$
or $\beta$ belongs to $(0,\frac12]$, then the kernels satisfy the
set of conditions (H1').

\begin{obs}
If we consider the kernel
processes $\theta_n^1$ and $\theta_n^2$, with the same L\'evy
process, it can be proved that we will obtain two family of approximation processes that
converge to a two independent fractional Brownian sheets. The proof
follows the  ideas in \cite{BMQ2020} for the case of the
Brownian sheet.
\end{obs}

\subsection{Other examples} In \cite{BJT2003} we can find other examples of
kernels that satisfy the set of conditions (H1').  So Theorem \ref{elteo}
 can be applied to processes with the
representation (\ref{fbs}) where $K_1$ and $K_2$ are some of these
kernels. For the sake of completeness let us recall these examples.
\subsubsection{Goursat kernels}
The kernel
$$K(t,r)=\sum_{i=1}^I g_i(t)h_i(r)I_{[0,t]}(r)$$ for some
$I\in\N$, with $g_i\in\textrm{Lip}_{\gamma_1} (0<\gamma_1\leq 1)$
and $h_i\in L^2([0,1])$. We impose also that $F$, defined by
$F(t)=\int_0^th_i^2(r)dr$, belongs to $\textrm{Lip}_{\gamma_2}
(0<\gamma_2\leq 1)$.

\subsubsection{The Holmgren-Riemann-Liouville fractional integral}
The kernel 
$$K(t,r)=\sqrt{2\pi}(t-r)^{H-\frac12}I_{[0,t]}(r),$$
with $0<H<1$. This kernel satisfy the set of conditions (H1') for
all $0<H<1$.
\subsubsection{A Lipschitz function}
The
kernel
$$K(t,r)=h(t-r)I_{[0,t]}(r),$$
with $h$ a Lipschitz function of parameter $\gamma\in(0,1]$.

\end{document}